\documentclass[a4paper]{article}

\usepackage{amssymb}

\newtheorem{thm}{Theorem}

\newtheorem{problem}{Problem}

\title{Bipartite graphs related to mutually  disjoint S-permutation matrices}
\author{Krasimir Yordzhev}
\date{}

\begin{document}
\unitlength=0.7mm \linethickness{0.5pt}
\maketitle

\begin{center} {\em
Faculty of Mathematics and Natural Sciences\\
South-West University, Blagoevgrad, Bulgaria} \\
E-mail: yordzhev@swu.bg
\end{center}

\begin{abstract}
Some numerical characteristics of bipartite graphs in relation to the problem of finding all disjoint pairs of S-permutation matrices in the general $n^2 \times n^2$ case are discussed in this paper. All  bipartite graphs of the type $g=\langle R_g \cup C_g , E_g \rangle$, where $|R_g |=|C_g |=2$ or $|R_g |=|C_g |=3$ are provided. The cardinality of the sets of mutually disjoint S-permutation matrices in both the $4 \times 4$ and $9 \times 9$ cases are calculated.
\end{abstract}

Keyword:{\it
 Bipartite graph, Binary matrix, S-permutation matrix, Disjoint matrices, Sudoku}

MSC[2010] code: 05C30, 05B20, 05C50

\section{Introduction}
Let $m$ be a positive integer. By $[m]$ we denote the set
$$[m] =\left\{ 1,2,\ldots ,m\right\} .$$

We let $\mathcal{S}_m$  denote  the symmetric group of order $m$ i.e., the group  of all one-to-one mappings of the  set $[m] $ to itself. If $x\in [m]$, $\rho\in \mathcal{S}_m$, then the image of the element $x$ in the mapping $\rho$ we will denote by $\rho (x)$.

A \emph{bipartite graph} is an ordered triple
 $$
g=\langle R_g , C_g , E_g \rangle ,
$$
where $R_g$ and $C_g$ are non-empty sets such that $R_g \cap C_g =\emptyset$. The elements of $R_g \cup C_g$ will be called \emph{vertices}. The set of \emph{edges} is $E_g \subseteq R_g \times C_g =\{ \langle r,c \rangle \; |\; r\in R_g ,c\in C_g \}$. Multiple edges are not allowed in our considerations.

The subject of the present work is bipartite graphs considered up to isomorphism.

 We refer to \cite{diestel} or  \cite{harary} for more details on graph theory.

Let $n$ and $k$ be two nonnegative  integers and let $0\le k\le n^2$. We denote by $\mathfrak{G}_{n,k} $ the set of all  bipartite graphs of the type $g=\langle  R_g , C_g ,E_g \rangle $, considered up to isomorphism, such that $ |R_g |= |C_g |= n$ and $|E_g |=k$.

Let $P_{ij}$, $1 \leq i,j \leq n$, be $n^2$ square $n\times n$ matrices, whose entries are elements of the set $[n^2 ] =\{ 1,2,\ldots ,n^2 \}$. The $n^2 \times n^2$ matrix

$$
P =
\left[
\begin{array}{cccc}
P_{11} & P_{12} & \cdots & P_{1n} \\
P_{21} & P_{22} & \cdots & P_{2n} \\
\vdots & \vdots & \ddots & \vdots \\
P_{n1} & P_{n2} & \cdots & P_{nn}
\end{array}
\right]
$$
is called a \emph{Sudoku matrix}, if every row, every column and every submatrix $P_{ij}$, $1\le i,j\le n$ comprise a permutation of the elements of set $[n^2 ]$, i.e., every number $s\in \{ 1,2,\ldots ,n^2 \}$ is found just once in each row, column, and submatrix $P_{ij}$. Submatrices $P_{ij}$ are called \emph{blocks} of $P$.

Sudoku is a very popular game and Sudoku matrices are special cases of Latin squares in the class of gerechte designs \cite{Bailey}.

A matrix is called \emph{binary} if all of its elements are equal to 0 or 1. A square binary matrix is called \emph{permutation matrix}, if in every row and every column there is just one 1.

Let us denote by $\Sigma_{n^2}$ the set of all $n^2 \times n^2$ permutation matrices of the following type:
$$
A =
\left[
\begin{array}{cccc}
A_{11} & A_{12} & \cdots & A_{1n} \\
A_{21} & A_{22} & \cdots & A_{2n} \\
\vdots & \vdots & \ddots & \vdots \\
A_{n1} & A_{n2} & \cdots & A_{nn}
\end{array}
\right] ,
$$
where for every $s,t\in \{ 1,2,\ldots ,n\}$,  $A_{st}$ is a  square $n\times n$ binary submatrix (block) with only one element equal to 1.

The elements of $\Sigma_{n^2} $ will be called \emph{S-permutation matrices}.

Two $\Sigma_{n^2}$  matrices $A=(a_{ij} )$ and $B=( b_{ij} )$, $1\le i,j\le n^2$ will be called \emph{disjoint}, if there are not elements $a_{ij}$ and $b_{ij}$ with the same indices such that $a_{ij} =b_{ij} =1$.

The concept of S-permutation matrix was introduced by Geir Dahl  \cite{dahl} in relation to the popular Sudoku puzzle.

Obviously, a square $n^2 \times n^2$ matrix $P$ with entries from $[n^2 ] =\{ 1,2,\ldots ,n^2 \}$ is a Sudoku matrix if and only if there are $\Sigma_{n^2}$ matrices $A_1 ,A_2 ,\ldots ,A_{n^2}$, pairwise disjoint, such that $P$ can be written in the following way:
\begin{equation}\label{disj}
P=1\cdot A_1 +2\cdot A_2 +\cdots +n^2 \cdot A_{n^2}
\end{equation}

In \cite{Fontana} Roberto Fontana offers an algorithm which returns a random family of $n^2 \times n^2$ mutually disjoint S-permutation matrices, where $n=2,3$. For $n=3$, he ran the algorithm 1000 times and found 105 different families of nine mutually disjoint S-permutation matrices. Then, applying  (\ref{disj}), he decided that there are at least $9!\cdot 105  =38\; 102\; 400$ Sudoku matrices. This number is very small compared with the exact number of $9\times 9$ Sudoku matrices. In  \cite{Felgenhauer} it was shown that  there are exactly
$$9! \cdot 72^2 \cdot 2^7 \cdot 27\; 704\; 267\; 971 = 6\; 670\; 903\; 752\; 021\; 072\; 936\; 960 $$
number of $9\times 9$ Sudoku matrices.

To evaluate the effectiveness of Fontana's algorithm, it is necessary to calculate the probability of two randomly generated matrices being disjoint.
As is proved in \cite{dahl}, the number of S-permutation matrices is equal to
$$
\left| \Sigma_{n^2} \right| = \left( n! \right)^{2n} .
$$

Thus the question of finding a  formula  for counting  disjoint pairs  of   S-permutation matrices naturally arises. Such a formula is introduced and verified in \cite{KYbigraphs}. In this paper, we demonstrate this formula  to  compute  the  number  of disjoint   pairs  of S-permutation matrices  in both  the 4 $\times$ 4 and 9 $\times$ 9 cases.

\section{A  formula  for counting  disjoint pairs  of   S-permutation matrices}
Let $g=\langle R_g , C_g ,E_g \rangle \in \mathfrak{G}_{n,k}$ for some natural numbers $n$ and $k$ and let $v\in V_g =R_g \cup C_g$.

By $N (v)$ we denote the set of all vertices of $V_g$, adjacent with $v$, i.e., $u\in N (v)$ if and only if there is an edge in $E_g$ connecting $u$ and $v$. If $v$ is an isolated vertex  (i.e., there is no edge, incident with $v$), then by definition $N (v)=\emptyset$ and $ \textrm{degree}  (v) = |N(v)|=0$. If $v\in R_g$, then obviously $N (v)\subseteq C_g$, and if $v\in C_g$, then $N (v)\subseteq R_g$.

Let $g=\langle R_g , C_g ,E_g \rangle \in \mathfrak{G}_{n,k}$ and let $u,v\in V_g =R_g\cup C_g$. We will say that $u$ and $v$ are equivalent and we will write $u\sim v$ if $N (u) =N (v)$. If $u$ and $v$ are isolated, then by definition $u\sim v$ if and only if $u$ and $v$ belong simultaneously to $R_g$, or $C_g$. The above introduced relation is obviously an equivalence relation.

By ${V_g}_{/\sim} $ we denote the obtained factor-set (the set of the equivalence classes) according to relation $\sim$ and let
$${V_g}_{/\sim} =\left\{ \Delta_1 ,\Delta_2 ,\ldots ,\Delta_s \right\} ,$$
where $\Delta_i \subseteq R_g$, or $\Delta_i \subseteq C_g$, $i=1,2,\ldots s$, $2\le s \le 2n$. We put
$$\delta_i =|\Delta_i |,\quad 1\le \delta_i \le n , \quad i=1,2,\ldots , s$$
and for every $g\in \mathfrak{G}_{n,k}$ we define multi-set (set with repetition)
$$\left[ g \right] =\left\{ \delta_1 ,\delta_2 ,\ldots \delta_s \right\} ,$$
where $\delta_1 ,\delta_2 ,\ldots ,\delta_s$ are natural numbers, obtained by the above described way.

If $z_1 \; z_2 \; \ldots \; z_n$ is a permutation of the elements of the set $[n] =\left\{ 1,2,\ldots ,n\right\}$ and we shortly denote $\rho$ this permutation, then in this case we denote by $\rho (i )$ the $i$-th element of this permutation, i.e., $\rho (i) =z_i$, $i=1,2,\ldots ,n$.

The following theorem is proved in \cite{KYbigraphs}:

\begin{thm}\label{t1} \cite{KYbigraphs}
Let $n\ge 2$ be a positive integer. Then the number $D_{n^2}$ of all disjoint ordered pairs of  matrices in $\Sigma_{n^2}$  is equal to
\begin{equation}\label{f2}
D_{n^2} = (n!)^{4n} + (n!)^{2(n+1)}\sum_{k=1}^{n^2} (-1)^k   \sum_{g\in \mathfrak{G}_{n,k} } \frac{\displaystyle \prod_{v\in R_g \cup C_g} ( n-|N (v)|)!}{\displaystyle \prod_{\delta \in [g]} \delta !}.
\end{equation}

The number  $d_{n^2}$ of all non-ordered pairs of disjoint matrices in $\Sigma_{n^2}$ is equal to
\begin{equation}\label{f3}
d_{n^2} =\frac{1}{2} D_{n^2}
\end{equation}

\hfill $\Box$
\end{thm}

The proof of Theorem \ref{t1} is described in detail in \cite{KYbigraphs} and here we will miss it.

In order to apply Theorem \ref{t1} it is necessary to describe all bipartite graphs up to isomorphism $g=\langle R_g , C_g , E_g \rangle $, where $|R_g |=|C_g |=n$.

Let $n$ and $k$ are positive integers and let $g\in \mathfrak{G}_{n,k}$. We examine the ordered $(n+1)$-tuple

\begin{equation}\label{PSI}
\Psi (g)=\langle \psi_0 (g) ,\psi_1 (g),\ldots ,\psi_n (g)\rangle ,
\end{equation}
where $\psi_i (g)$, $i=0,1,\ldots ,n$ is equal to the number of vertices of $g$ incident with exactly $i$ number of edges. It is obvious that $\displaystyle \sum_{i=1}^n i\psi_i (g)=2k$ is true for all $g\in \mathfrak{G}_{n,k}$. Then formula (\ref{f2}) can be presented
$$
D_{n^2} = (n!)^{4n} + (n!)^{2(n+1)}\sum_{k=1}^{n^2} (-1)^k   \sum_{g\in \mathfrak{G}_{n,k} } \frac{\displaystyle \prod_{i=0}^n \left[ \left( n-i\right) ! \right]^{\psi_i (g)}}{\displaystyle \prod_{\delta \in [g]} \delta !}.
$$

Since $(n-n)!=0!=1$ and $[n-(n-1)]!=1!=1 $, then

\begin{equation}\label{f2star}
D_{n^2} = (n!)^{4n} + (n!)^{2(n+1)}\sum_{k=1}^{n^2} (-1)^k   \sum_{g\in \mathfrak{G}_{n,k} } \frac{\displaystyle \prod_{i=0}^{n-2} \left[ \left( n-i\right) ! \right]^{\psi_i (g)}}{\displaystyle \prod_{\delta \in [g]} \delta !}.
\end{equation}

Consequently, to apply formula (\ref{f2star}) for each bipartite graph $g\in \mathfrak{G}_{n,k}$ and for the set $\mathfrak{G}_{n,k}$ of  bipartite graphs, it is necessary to obtain the following numerical characteristics:

\begin{equation}\label{star2}
\omega (g) = \frac{\displaystyle \prod_{i=0}^{n-2} \left[ \left( n-i\right) ! \right]^{\psi_i (g)}}{\displaystyle \prod_{\delta \in [g]} \delta !}
\end{equation}
and
\begin{equation}\label{thet}
\theta (n,k) =\sum_{g\in \mathfrak{G}_{n,k} } \omega (g)
\end{equation}

Using the numerical characteristics (\ref{star2}) and (\ref{thet}), we obtain the following variety  of Theorem \ref{t1}:

\begin{thm}\label{t2}
\begin{equation}\label{endstar}
D_{n^2} = (n!)^{4n} + (n!)^{2(n+1)}\sum_{k=1}^{n^2} (-1)^k   \theta (n,k),
\end{equation}
where $\theta (n,k)$ is described using formulas (\ref{thet}) and (\ref{star2}).

\hfill $\Box$
\end{thm}

\section{Demonstrations in applying of Theorem   \ref{t2}}

\subsection{Counting the number $D_4$ of all ordered pairs of disjoint S-permutation matrices for $n = 2$}

\subsubsection{$k=1$}
In $n=2$ and $k=1$, $\mathfrak{G}_{2,1}$ consists of a single graph $g_1$ shown in Figure \ref{n2k1}.

\begin{figure}[h]
\begin{picture}(50,30)

\put(0,0.01){\framebox(50,30)}
\put(3,27){\makebox(0,0){$g_1$}}

\multiput(13,16)(24,0){2}{\oval(13,20)}
\put(14,3){\makebox(0,0){$R_{g_1}$}}
\put(38,3){\makebox(0,0){$C_{g_1}$}}

\put(13,20){\circle*{2}}
\put(37,20){\circle{2}}
\put(15,20){\line(1,0){20}}

\put(13,10){\circle*{2}}
\put(37,10){\circle{2}}

\end{picture}
\caption{$n=2$, $k=1$}\label{n2k1}
\end{figure}

For graph $g_1 \in \mathfrak{G}_{2,1}$ we have:
$$\left[ g_1 \right] =\left\{ 1,1,1,1\right\}$$
$$\Psi (g_1 )=\langle \psi_0 (g_1 ) ,\psi_1 (g_1 ),\psi_2 (g_1 )\rangle =\langle 2,2,0\rangle$$

Then we get:
$$\omega (g_1 )=\frac{\left[ (2-0)!\right]^2 }{1!\; 1!\; 1!\; 1!} =4$$
and therefore
\begin{equation}\label{theta21}
\theta (2,1) =\sum_{g\in\mathfrak{G}_{2,1}} \omega (g) =4 .
\end{equation}

\subsubsection{$k=2$}
The set $\mathfrak{G}_{2,2}$ consists of three graphs $g_2$, $g_3$ and $g_4$ depicted in Figure \ref{n2k2}.

\begin{figure}[h]
\begin{picture}(170,30)

\multiput(0,0)(60,0){3}{
\put(0,0.01){\framebox(50,30)}

\multiput(13,16)(24,0){2}{\oval(13,20)}

\put(13,20){\circle*{2}}
\put(37,20){\circle{2}}

\put(13,10){\circle*{2}}
\put(37,10){\circle{2}}
}

\put(15,20){\line(1,0){20}}
\put(15,10){\line(1,0){20}}

\put(75,20){\line(1,0){20}}
\put(75,10){\line(2,1){20}}

\put(135,20){\line(1,0){20}}
\put(135,20){\line(2,-1){20}}

\put(3,27){\makebox(0,0){$g_2$}}
\put(63,27){\makebox(0,0){$g_3$}}
\put(123,27){\makebox(0,0){$g_4$}}

\put(14,3){\makebox(0,0){$R_{g_2}$}}
\put(38,3){\makebox(0,0){$C_{g_2}$}}

\put(74,3){\makebox(0,0){$R_{g_3}$}}
\put(98,3){\makebox(0,0){$C_{g_3}$}}

\put(134,3){\makebox(0,0){$R_{g_4}$}}
\put(158,3){\makebox(0,0){$C_{g_4}$}}

\end{picture}
\caption{$n=2$, $k=2$}\label{n2k2}
\end{figure}

For graph $g_2 \in \mathfrak{G}_{2,2}$ we have:
$$\left[ g_2 \right] =\left\{ 1,1,1,1\right\}$$
$$\Psi (g_2 )=\langle \psi_0 (g_2 ) ,\psi_1 (g_2 ),\psi_2 (g_2 )\rangle =\langle 0,4,0\rangle$$

$$\omega (g_1 )=\frac{\left[ (2-0)!\right]^0 }{1!\; 1!\; 1!\; 1!} =1$$

For graphs $g_3 \in \mathfrak{G}_{2,2}$ and $g_4 \in \mathfrak{G}_{2,2}$ we have:
$$\left[ g_3 \right] = \left[ g_4 \right] =\left\{ 2,1,1\right\}$$
$$\Psi (g_3 ) =\Psi (g_4 ) =\langle 1,2,1\rangle$$
$$\omega (g_3 )=\omega (g_4 )=\frac{\left[ (2-0)!\right]^1 }{2!\; 1!\; 1!}  =1$$

Then for the set $\mathfrak{G}_{2,2}$ we get:
\begin{equation}\label{theta22}
\theta (2,2) =\sum_{g\in\mathfrak{G}_{2,2}} \omega (g) = 1+1+1=3 .
\end{equation}

\subsubsection{$k=3$}
In $n = 2$ and $k = 3$, $\mathfrak{G}_{2,3}$ consists of a single graph $g_5$ shown in Figure \ref{n2k3}.

\begin{figure}[h]
\begin{picture}(50,30)

\put(0,0.01){\framebox(50,30)}
\put(3,27){\makebox(0,0){$g_5$}}

\multiput(13,16)(24,0){2}{\oval(13,20)}
\put(14,3){\makebox(0,0){$R_{g_5}$}}
\put(38,3){\makebox(0,0){$C_{g_5}$}}

\put(13,20){\circle*{2}}
\put(37,20){\circle{2}}

\put(15,20){\line(2,-1){20}}

\put(13,10){\circle*{2}}
\put(37,10){\circle{2}}
\put(15,10){\line(1,0){20}}
\put(15,10){\line(2,1){20}}

\end{picture}
\caption{$n=2$, $k=3$}\label{n2k3}
\end{figure}

For graph $g_5 \in \mathfrak{G}_{2,3}$ we have:
$$\left[ g_5 \right] =\left\{ 1,1,1,1\right\}$$
$$\Psi (g_5 )=\langle \psi_0 (g_5 ) ,\psi_1 (g_5 ),\psi_2 (g_5 )\rangle =\langle 0,2,2\rangle$$

Then we get:
$$\omega (g_5 )=\frac{\left[ (2-0)!\right]^0 }{1!\; 1!\; 1!\; 1!}  =1$$
and therefore
\begin{equation}\label{theta23}
\theta (2,3) =\sum_{g\in\mathfrak{G}_{2,3}} \omega (g) =1 .
\end{equation}

\subsubsection{$k=4$}

When $n = 2$ and $k = 4$ there is only one graph and this is the complete bipartite graph $g_6$ which is shown in Figure \ref{n2k4}.

\begin{figure}[h]
\begin{picture}(50,30)

\put(0,0.01){\framebox(50,30)}
\put(3,27){\makebox(0,0){$g_6$}}

\multiput(13,16)(24,0){2}{\oval(13.5,20)}
\put(14,3){\makebox(0,0){$R_{g_6}$}}
\put(38,3){\makebox(0,0){$C_{g_6}$}}

\put(13,20){\circle*{2}}
\put(37,20){\circle{2}}
\put(15,20){\line(1,0){20}}
\put(15,20){\line(2,-1){20}}

\put(13,10){\circle*{2}}
\put(37,10){\circle{2}}
\put(15,10){\line(1,0){20}}
\put(15,10){\line(2,1){20}}

\end{picture}
\caption{$n=2$, $k=4$}\label{n2k4}
\end{figure}

For graph $g_6 \in \mathfrak{G}_{2,4}$ we have:

$$\left[ g_6 \right] =\left\{ 2,2\right\}$$
$$\Psi (g_6 )=\langle \psi_0 (g_6 ) ,\psi_1 (g_6 ),\psi_2 (g_6 )\rangle =\langle 0,0,4\rangle$$

Then we get:
$$\omega (g_6 )=\frac{\left[ (2-0)!\right]^0  }{2!\; 2!}  =\frac{1}{4}$$
and therefore
\begin{equation}\label{theta24}
\theta (2,4) =\sum_{g\in\mathfrak{G}_{2,1}} \omega (g) =\frac{1}{4} .
\end{equation}

Having in mind the formulas (\ref{endstar}), (\ref{theta21}), (\ref{theta22}), (\ref{theta23}) and (\ref{theta24}) for the number $D_4$ of all ordered pairs disjoint S-permutation  matrices in $ n =  2$ we finally  get:

\begin{equation}\label{D_4}
D_4 = (2!)^{8} + (2!)^{6} \left[ -\theta (2,1)+\theta (2,2)-\theta (2,3)+\theta (2,4) \right] =
\end{equation}
$$=256+64 \left( -4+3-1+\frac{1}{4} \right) =144.$$

The number $d_4$ of all non-ordered pairs disjoint matrices from $\Sigma_4$  is equal to

\begin{equation}\label{ddd2}
d_4 =\frac{1}{2} D_4 =72.
\end{equation}

\subsection{Counting the number $D_9$ of all ordered pairs of disjoint S-permutation matrices for $n = 3$ }

\subsubsection{$k=1$}

Graph $g_7$, which is displayed in Figure \ref{n3k1} is the only bipartite graph belonging to the set $\mathfrak{G}_{3,1}$.

\begin{figure}[h]
\begin{picture}(50,40)

\put(0,0){\framebox(50,40)}
\put(3,37){\makebox(0,0){$g_7$}}

\multiput(13,21)(24,0){2}{\oval(14,30)}
\put(14,3){\makebox(0,0){$R_{g_7}$}}
\put(38,3){\makebox(0,0){$C_{g_7}$}}

\put(13,30){\circle*{2}}
\put(37,30){\circle{2}}
\put(15,30){\line(1,0){20}}

\put(13,20){\circle*{2}}
\put(37,20){\circle{2}}

\put(13,10){\circle*{2}}
\put(37,10){\circle{2}}

\end{picture}
\caption{$n=3$, $k=1$}\label{n3k1}
\end{figure}

For graph $g_7 \in \mathfrak{G}_{3,1}$ we have:

$$\left[ g_7 \right] =\left\{ 1,1,2,2\right\}$$
$$\Psi (g_7 )=\langle \psi_0 (g_7 ) ,\psi_1 (g_7 ),\psi_2 (g_7 ),\psi_3 (g_7 ),\psi_4 (g_8 )\rangle =\langle 4,2,0,0\rangle$$

Then we get:

$$\omega (g_7 )=\frac{\left[ (3-0)!\right]^4 \left[ (3-1)!\right]^2 }{1!\; 1!\; 2!\; 2!} =\frac{6^4\cdot 2^2 }{1\cdot 1\cdot 2\cdot 2} =1296$$
and therefore
\begin{equation}\label{theta31}
\theta (3,1) =\sum_{g\in\mathfrak{G}_{3,1}} \omega (g) =1296 .
\end{equation}

\subsubsection{$k=2$}

In this case $\mathfrak{G}_{3,2} =\{ g_8 ,g_9 ,g_{10} \}$. The graphs $ g_8 $, $ g_9 $ and $ g_{10} $ are shown in Figure \ref{n3k2}.

\begin{figure}[h]
\begin{picture}(170,40)

\multiput(0,0)(60,0){3}{
\put(0,0){\framebox(50,40)}

\multiput(13,21)(24,0){2}{\oval(14,30)}

\put(13,30){\circle*{2}}
\put(37,30){\circle{2}}

\put(13,20){\circle*{2}}
\put(37,20){\circle{2}}

\put(13,10){\circle*{2}}
\put(37,10){\circle{2}}
}

\put(3,37){\makebox(0,0){$g_8$}}
\put(63,37){\makebox(0,0){$g_9$}}
\put(123,37){\makebox(0,0){$g_{10}$}}

\put(14,3){\makebox(0,0){$R_{g_8}$}}
\put(38,3){\makebox(0,0){$C_{g_8}$}}
\put(74,3){\makebox(0,0){$R_{g_9}$}}
\put(98,3){\makebox(0,0){$C_{g_9}$}}
\put(134,3){\makebox(0,0){$R_{g_{10}}$}}
\put(158,3){\makebox(0,0){$C_{g_{10}}$}}

\put(15,30){\line(1,0){20}}
\put(15,20){\line(1,0){20}}

\put(75,30){\line(1,0){20}}
\put(75,30){\line(2,-1){20}}

\put(135,30){\line(1,0){20}}
\put(135,20){\line(2,1){20}}

\end{picture}
\caption{$n=3$, $k=2$}\label{n3k2}
\end{figure}

For graph $g_8 \in \mathfrak{G}_{3,2}$ we have:

$$\left[ g_8 \right] =\left\{ 1,1,1,1,1,1\right\}$$

$$\Psi (g_8 )=\langle \psi_0 (g_8 ) ,\psi_1 (g_8 ),\psi_2 (g_8 ),\psi_3 (g_8 ),\psi_4 (g_8 )\rangle =\langle 2,4,0,0\rangle$$

$$\omega (g_8 )=\frac{\left[ (3-0)!\right]^2 \left[ (3-1)!\right]^4 }{1!\; 1!\; 1!\; 1!\; 1!\; 1!} = 6^2\cdot 2^4 =576$$

For graphs $g_9 \in \mathfrak{G}_{3,2}$ and $g_{10} \in \mathfrak{G}_{3,2}$ we have:

$$\left[ g_9 \right] =[g_{10} ]=\left\{ 1,1,2,2\right\}$$

$$\Psi (g_9 )= \Psi (g_{10} ) =\langle 3,2,1,0\rangle$$

$$\omega (g_9 )=\omega (g_{10} )=\frac{\left[ (3-0)!\right]^3 \left[ (3-1)!\right]^2 }{1!\; 1!\; 2!\; 2!} = \frac{6^3 \cdot 2^2}{1\cdot 1\cdot 2\cdot 2} =216$$

Then for the set $\mathfrak{G}_{3,2}$ we get:

\begin{equation}\label{theta32}
\theta (3,2) =\sum_{g\in\mathfrak{G}_{3,2}} \omega (g) =576+216+216=1008 .
\end{equation}

\subsubsection{$k=3$}

When $n = 3$ and $k = 3$ the set $\mathfrak{G}_{3,3}= \{ g_{11} ,g_{12} ,g_{13} ,g_{14} ,g_{15} ,g_{16} \}$ consists of six bipartite graphs, which are shown in Figure \ref{n3k3}.

\begin{figure}[h]
\begin{picture}(170,90)

\multiput(0,0)(60,0){3}{
\multiput(0,0)(0,50){2}{

\put(0,0){\framebox(50,40)}

\multiput(13,21)(24,0){2}{\oval(14,30)}

\put(13,30){\circle*{2}}
\put(37,30){\circle{2}}

\put(13,20){\circle*{2}}
\put(37,20){\circle{2}}

\put(13,10){\circle*{2}}
\put(37,10){\circle{2}}
}
}

\put(3,87){\makebox(0,0){$g_{11}$}}
\put(63,87){\makebox(0,0){$g_{12}$}}
\put(123,87){\makebox(0,0){$g_{13}$}}

\put(14,53){\makebox(0,0){$R_{g_{11}}$}}
\put(38,53){\makebox(0,0){$C_{g_{11}}$}}

\put(74,53){\makebox(0,0){$R_{g_{12}}$}}
\put(98,53){\makebox(0,0){$C_{g_{12}}$}}

\put(134,53){\makebox(0,0){$R_{g_{13}}$}}
\put(158,53){\makebox(0,0){$C_{g_{13}}$}}

\put(3,37){\makebox(0,0){$g_{14}$}}
\put(63,37){\makebox(0,0){$g_{15}$}}
\put(123,37){\makebox(0,0){$g_{16}$}}

\put(14,3){\makebox(0,0){$R_{g_{14}}$}}
\put(38,3){\makebox(0,0){$C_{g_{14}}$}}

\put(74,3){\makebox(0,0){$R_{g_{15}}$}}
\put(98,3){\makebox(0,0){$C_{g_{15}}$}}

\put(134,3){\makebox(0,0){$R_{g_{16}}$}}
\put(158,3){\makebox(0,0){$C_{g_{16}}$}}

\put(15,80){\line(1,0){20}}
\put(15,70){\line(1,0){20}}
\put(15,60){\line(1,0){20}}

\put(75,80){\line(1,0){20}}
\put(75,80){\line(2,-1){20}}
\put(75,60){\line(1,0){20}}

\put(135,80){\line(1,0){20}}
\put(135,70){\line(2,1){20}}
\put(135,60){\line(1,0){20}}

\put(15,30){\line(1,0){20}}
\put(15,30){\line(2,-1){20}}
\put(15,20){\line(1,0){20}}

\put(75,20){\line(1,0){20}}
\put(75,20){\line(2,1){20}}
\put(75,20){\line(2,-1){20}}

\put(135,30){\line(2,-1){20}}

\put(135,20){\line(1,0){20}}

\put(135,10){\line(2,1){20}}

\end{picture}
\caption{$n=3$, $k=3$}\label{n3k3}
\end{figure}

For graph $g_{11} \in\mathfrak{G}_{3,3}$ we have:

$$\left[ g_{11} \right] =\left\{ 1,1,1,1,1,1\right\}$$

$$\Psi (g_{11} ) =\langle 0,6,0,0\rangle$$

$$\omega (g_{11} )=\frac{\left[ (3-0)!\right]^0 \left[ (3-1)!\right]^6 }{1!\; 1!\; 1!\; 1!\; 1!\; 1!} = 6^0 \cdot 2^6 =64$$

For graphs $g_{12} ,g_{13} \in\mathfrak{G}_{3,3}$ we have:

$$\left[ g_{12} \right] =\left[ g_{13} \right] =\left\{ 1,1,1,1,2\right\}$$

$$\Psi (g_{12} ) =\Psi (g_{13} ) =\langle 1,4,1,0\rangle$$

$$\omega (g_{12} )=\omega (g_{13} )=\frac{\left[ (3-0)!\right]^1 \left[ (3-1)!\right]^4 }{1!\; 1!\; 1!\; 1!\; 2!} = \frac{6^1 \cdot 2^4}{2} =48$$

For graph $g_{14} \in\mathfrak{G}_{3,3}$ we have:

$$\left[ g_{14} \right] =\left\{ 1,1,1,1,1,1\right\}$$

$$\Psi (g_{14} ) =\langle 2,2,2,0\rangle$$

$$\omega (g_{14} )=\frac{\left[ (3-0)!\right]^2 \left[ (3-1)!\right]^2 }{1!\; 1!\; 1!\; 1!\; 1!\; 1!} = 6^2 \cdot 2^2 =144$$

For graphs $g_{15} ,g_{16} \in\mathfrak{G}_{3,3}$ we have:

$$\left[ g_{15} \right] =\left[ g_{16} \right] =\left\{ 1,2,3\right\}$$

$$\Psi (g_{15} ) =\Psi (g_{16} ) =\langle 2,3,0,1\rangle$$

$$\omega (g_{15} )=\omega (g_{16} )=\frac{\left[ (3-0)!\right]^2 \left[ (3-1)!\right]^3 }{1!\; 2!\; 3!} = \frac{6^2 \cdot 2^3}{2\cdot 6} =24$$

Then for the set $\mathfrak{G}_{3,3}$ we get:

\begin{equation}\label{theta33}
\theta (3,3)=\sum_{g\in \mathfrak{G}_{3,3}} \omega (g)=64+48+48+144+24+24=352 .
\end{equation}

\subsubsection{$k=4$}

When $n = 3$ and $k = 4$ the set $\mathfrak{G}_{3,4}= \{ g_{17} ,g_{18} ,g_{19} ,g_{20} ,g_{21} ,g_{22} ,g_{23} \}$ consists of seven bipartite graphs, which are shown in Figure \ref{n3k4}.

\begin{figure}[h]
\begin{picture}(170,140)

\multiput(0,0)(60,0){2}{
\multiput(0,0)(0,50){3}{

\put(0,0){\framebox(50,40)}

\multiput(13,21)(24,0){2}{\oval(14,30)}

\put(13,30){\circle*{2}}
\put(37,30){\circle{2}}

\put(13,20){\circle*{2}}
\put(37,20){\circle{2}}

\put(13,10){\circle*{2}}
\put(37,10){\circle{2}}
}
}
\put(120,100){\framebox(50,40)}

\multiput(133,121)(24,0){2}{\oval(14,30)}

\put(133,130){\circle*{2}}
\put(157,130){\circle{2}}

\put(133,120){\circle*{2}}
\put(157,120){\circle{2}}

\put(133,110){\circle*{2}}
\put(157,110){\circle{2}}

\put(3,137){\makebox(0,0){$g_{17}$}}
\put(63,137){\makebox(0,0){$g_{18}$}}
\put(123,137){\makebox(0,0){$g_{19}$}}

\put(14,103){\makebox(0,0){$R_{g_{17}}$}}
\put(38,103){\makebox(0,0){$C_{g_{17}}$}}

\put(74,103){\makebox(0,0){$R_{g_{18}}$}}
\put(98,103){\makebox(0,0){$C_{g_{18}}$}}

\put(134,103){\makebox(0,0){$R_{g_{19}}$}}
\put(158,103){\makebox(0,0){$C_{g_{19}}$}}

\put(3,87){\makebox(0,0){$g_{20}$}}
\put(63,87){\makebox(0,0){$g_{21}$}}

\put(14,53){\makebox(0,0){$R_{g_{20}}$}}
\put(38,53){\makebox(0,0){$C_{g_{20}}$}}

\put(74,53){\makebox(0,0){$R_{g_{21}}$}}
\put(98,53){\makebox(0,0){$C_{g_{21}}$}}

\put(3,37){\makebox(0,0){$g_{22}$}}
\put(63,37){\makebox(0,0){$g_{23}$}}

\put(14,3){\makebox(0,0){$R_{g_{22}}$}}
\put(38,3){\makebox(0,0){$C_{g_{22}}$}}

\put(74,3){\makebox(0,0){$R_{g_{23}}$}}
\put(98,3){\makebox(0,0){$C_{g_{23}}$}}

\put(15,130){\line(1,0){20}}
\put(15,130){\line(2,-1){20}}

\put(15,120){\line(1,0){20}}
\put(15,120){\line(2,1){20}}


\put(75,130){\line(1,0){20}}
\put(75,130){\line(2,-1){20}}

\put(75,120){\line(2,-1){20}}

\put(75,110){\line(1,0){20}}

\put(135,130){\line(1,0){20}}
\put(135,130){\line(2,-1){20}}

\put(135,120){\line(1,0){20}}

\put(135,110){\line(1,0){20}}

\put(15,80){\line(1,0){20}}

\put(15,70){\line(2,1){20}}
\put(15,70){\line(2,-1){20}}

\put(15,60){\line(1,0){20}}

\put(75,80){\line(1,0){20}}
\put(75,80){\line(2,-1){20}}


\put(75,60){\line(1,0){20}}
\put(75,60){\line(2,1){20}}

\put(15,30){\line(1,0){20}}

\put(15,20){\line(1,0){20}}
\put(15,20){\line(2,1){20}}
\put(15,20){\line(2,-1){20}}


\put(75,30){\line(1,0){20}}
\put(75,30){\line(2,-1){20}}

\put(75,20){\line(1,0){20}}

\put(75,10){\line(2,1){20}}

\end{picture}
\caption{$n=3$, $k=4$}\label{n3k4}
\end{figure}

For graph $g_{17} \in\mathfrak{G}_{3,4}$ we have:

$$\left[ g_{17} \right] =\left\{ 1,1,2,2\right\}$$

$$\Psi (g_{17} ) =\langle 2,0,4,0\rangle$$

$$\omega (g_{17} )=\frac{\left[ (3-0)!\right]^2 \left[ (3-1)!\right]^0 }{1!\; 1!\; 2!\; 2!} = \frac{6^2 \cdot 2^0}{2^2} =9$$

For graph $g_{18} \in\mathfrak{G}_{3,4}$ we have:

$$\left[ g_{18} \right] =\left\{ 1,1,2,2\right\}$$

$$\Psi (g_{18} ) =\langle 0,4,2,0\rangle$$

$$\omega (g_{18} )=\frac{\left[ (3-0)!\right]^0 \left[ (3-1)!\right]^4 }{1!\; 1!\; 2!\; 2!} = \frac{6^0 \cdot 2^4}{2^2} =4$$

For graph $g_{19} \in\mathfrak{G}_{3,4}$ we have:

$$\left[ g_{19} \right] =\left\{ 1,1,1,1,1,1\right\}$$

$$\Psi (g_{19} ) =\langle 0,4,2,0\rangle$$

$$\omega (g_{19} )=\frac{\left[ (3-0)!\right]^0 \left[ (3-1)!\right]^4 }{1!\; 1!\; 1!\; 1!\; 1!\; 1!} = 6^0 \cdot 2^4 =16$$

For graphs $g_{20} \in \mathfrak{G}_{3,4}$ and $g_{21} \in \mathfrak{G}_{3,4}$ we have:

$$\left[ g_{20} \right] =\left[ g_{21} \right] =\left\{ 1,1,1,1,1,1\right\}$$

$$\Psi (g_{20} ) =\Psi (g_{21} ) =\langle 1,2,3,0\rangle$$

$$\omega (g_{20} )=\omega (g_{21} )=\frac{\left[ (3-0)!\right]^1 \left[ (3-1)!\right]^2 }{1!\; 1!\; 1!\; 1!\; 1!\; 1!} = 6^1 \cdot 2^2 =24$$

For graphs $g_{22} \in \mathfrak{G}_{3,4}$ and $g_{23} \in \mathfrak{G}_{3,4}$ we have:

$$\left[ g_{22} \right] =\left[ g_{23} \right] =\left\{ 1,1,1,1,2\right\}$$

$$\Psi (g_{22} ) =\Psi (g_{23} ) =\langle 1,3,1,1\rangle$$

$$\omega (g_{22} )=\omega (g_{23} )=\frac{\left[ (3-0)!\right]^1 \left[ (3-1)!\right]^3 }{1!\; 1!\; 1!\; 1!\; 2!} = \frac{6^1 \cdot 2^3}{2} =24$$

Then we get:

\begin{equation}\label{theta34}
\theta (3,4)=\sum_{g\in \mathfrak{G}_{3,4}} \omega (g)=9+4+16+24+24+24+24=125 .
\end{equation}

\subsubsection{$k=5$}

When $n = 3$ and $k = 5$ the set $\mathfrak{G}_{3,5}$ consists of seven bipartite graphs $g_{24} \div g_{30}$, which are shown in Figure  \ref{n3k5}.

\begin{figure}[h]
\begin{picture}(170,140)

\multiput(0,0)(60,0){2}{
\multiput(0,0)(0,50){3}{

\put(0,0){\framebox(50,40)}

\multiput(13,21)(24,0){2}{\oval(14,30)}

\put(13,30){\circle*{2}}
\put(37,30){\circle{2}}

\put(13,20){\circle*{2}}
\put(37,20){\circle{2}}

\put(13,10){\circle*{2}}
\put(37,10){\circle{2}}
}
}
\put(120,100){\framebox(50,40)}

\multiput(133,121)(24,0){2}{\oval(14,30)}

\put(133,130){\circle*{2}}
\put(157,130){\circle{2}}

\put(133,120){\circle*{2}}
\put(157,120){\circle{2}}

\put(133,110){\circle*{2}}
\put(157,110){\circle{2}}

\put(3,137){\makebox(0,0){$g_{24}$}}
\put(63,137){\makebox(0,0){$g_{25}$}}
\put(123,137){\makebox(0,0){$g_{26}$}}

\put(14,103){\makebox(0,0){$R_{g_{24}}$}}
\put(38,103){\makebox(0,0){$C_{g_{24}}$}}

\put(74,103){\makebox(0,0){$R_{g_{25}}$}}
\put(98,103){\makebox(0,0){$C_{g_{25}}$}}

\put(134,103){\makebox(0,0){$R_{g_{26}}$}}
\put(158,103){\makebox(0,0){$C_{g_{26}}$}}

\put(3,87){\makebox(0,0){$g_{27}$}}
\put(63,87){\makebox(0,0){$g_{28}$}}

\put(14,53){\makebox(0,0){$R_{g_{27}}$}}
\put(38,53){\makebox(0,0){$C_{g_{27}}$}}

\put(74,53){\makebox(0,0){$R_{g_{28}}$}}
\put(98,53){\makebox(0,0){$C_{g_{28}}$}}

\put(3,37){\makebox(0,0){$g_{29}$}}
\put(63,37){\makebox(0,0){$g_{30}$}}

\put(14,3){\makebox(0,0){$R_{g_{29}}$}}
\put(38,3){\makebox(0,0){$C_{g_{29}}$}}

\put(74,3){\makebox(0,0){$R_{g_{30}}$}}
\put(98,3){\makebox(0,0){$C_{g_{30}}$}}

\put(15,130){\line(1,-1){20}}

\put(15,120){\line(2,-1){20}}

\put(15,110){\line(1,0){20}}
\put(15,110){\line(1,1){20}}
\put(15,110){\line(2,1){20}}

\put(75,130){\line(1,-1){20}}

\put(75,120){\line(1,0){20}}
\put(75,120){\line(2,1){20}}

\put(75,110){\line(1,1){20}}
\put(75,110){\line(2,1){20}}

\put(135,130){\line(1,-1){20}}

\put(135,120){\line(2,1){20}}
\put(135,120){\line(2,-1){20}}

\put(135,110){\line(1,1){20}}
\put(135,110){\line(2,1){20}}

\put(15,80){\line(1,-1){20}}
\put(15,80){\line(2,-1){20}}

\put(15,70){\line(1,0){20}}

\put(15,60){\line(1,1){20}}
\put(15,60){\line(2,1){20}}

\put(75,80){\line(1,-1){20}}

\put(75,70){\line(1,0){20}}
\put(75,70){\line(2,1){20}}
\put(75,70){\line(2,-1){20}}

\put(75,60){\line(1,1){20}}

\put(15,30){\line(1,-1){20}}
\put(15,30){\line(2,-1){20}}


\put(15,10){\line(1,0){20}}
\put(15,10){\line(1,1){20}}
\put(15,10){\line(2,1){20}}

\put(75,30){\line(1,-1){20}}

\put(75,20){\line(2,1){20}}
\put(75,20){\line(2,-1){20}}

\put(75,10){\line(1,0){20}}
\put(75,10){\line(1,1){20}}

\end{picture}
\caption{$n=3$, $k=5$}\label{n3k5}
\end{figure}

For graph $g_{24} \in\mathfrak{G}_{3,5}$ we have:

$$\left[ g_{24} \right] =\left\{ 1,1,2,2\right\}$$

$$\Psi (g_{24} ) =\langle 0,4,0,2\rangle$$

$$\omega (g_{18} )=\frac{\left[ (3-0)!\right]^0 \left[ (3-1)!\right]^4 }{1!\; 1!\; 2!\; 2!} = \frac{6^0 \cdot 2^4}{2^2} =4$$

For graph $g_{25} \in\mathfrak{G}_{3,5}$ we have:

$$\left[ g_{25} \right] =\left\{ 1,1,2,2\right\}$$

$$\Psi (g_{25} ) =\langle 0,2,4,0\rangle$$

$$\omega (g_{18} )=\frac{\left[ (3-0)!\right]^0 \left[ (3-1)!\right]^2 }{1!\; 1!\; 2!\; 2!} = \frac{6^0 \cdot 2^2}{2^2} =1$$

For graph $g_{26} \in\mathfrak{G}_{3,5}$ we have:

$$\left[ g_{26} \right] =\left\{ 1,1,1,1,1,1\right\}$$

$$\Psi (g_{26} ) =\langle 0,2,4,0\rangle$$

$$\omega (g_{26} )=\frac{\left[ (3-0)!\right]^0 \left[ (3-1)!\right]^2 }{1!\; 1!\; 1!\; 1!\; 1!\; 1!} = 6^0 \cdot 2^2 =4$$

For graphs $g_{27} \in \mathfrak{G}_{3,5}$ and $g_{28} \in \mathfrak{G}_{3,5}$ we have:

$$\left[ g_{27} \right] =\left[ g_{28} \right] =\left\{ 1,1,1,1,1,1\right\}$$

$$\Psi (g_{27} ) =\Psi (g_{28} ) =\langle 0,3,2,1\rangle$$

$$\omega (g_{27} )=\omega (g_{28} )=\frac{\left[ (3-0)!\right]^0 \left[ (3-1)!\right]^3 }{1!\; 1!\; 1!\; 1!\; 1!\; 1!} = 6^0 \cdot 2^3 =8$$

For graphs $g_{29} \in \mathfrak{G}_{3,5}$ and $g_{30} \in \mathfrak{G}_{3,5}$ we have:

$$\left[ g_{29} \right] =\left[ g_{30} \right] =\left\{ 1,1,1,1,2\right\}$$

$$\Psi (g_{29} ) =\Psi (g_{30} ) =\langle 1,1,3,1\rangle$$

$$\omega (g_{29} )=\omega (g_{30} )=\frac{\left[ (3-0)!\right]^1 \left[ (3-1)!\right]^1 }{1!\; 1!\; 1!\; 1!\; 2!} = \frac{6^1 \cdot 2^1}{2} =6$$

Then we get:

\begin{equation}\label{theta35}
\theta (3,5)=\sum_{g\in \mathfrak{G}_{3,5}} \omega (g)=4+1+4+8+8+6+6=37 .
\end{equation}

\subsubsection{$k=6$}
When $n = 3$ and $k = 6$ the set $\mathfrak{G}_{3,6}= \{ g_{31} ,g_{32} ,g_{33} ,g_{34} ,g_{35} ,g_{36} \}$ consists of six bipartite graphs, which are shown in Figure \ref{n3k6}.

\begin{figure}[h]
\begin{picture}(170,90)

\multiput(0,0)(60,0){3}{
\multiput(0,0)(0,50){2}{

\put(0,0){\framebox(50,40)}

\multiput(13,21)(24,0){2}{\oval(14,30)}

\put(13,30){\circle*{2}}
\put(37,30){\circle{2}}

\put(13,20){\circle*{2}}
\put(37,20){\circle{2}}

\put(13,10){\circle*{2}}
\put(37,10){\circle{2}}
}
}

\put(3,87){\makebox(0,0){$g_{31}$}}
\put(63,87){\makebox(0,0){$g_{32}$}}
\put(123,87){\makebox(0,0){$g_{33}$}}

\put(14,53){\makebox(0,0){$R_{g_{31}}$}}
\put(38,53){\makebox(0,0){$C_{g_{31}}$}}

\put(74,53){\makebox(0,0){$R_{g_{32}}$}}
\put(98,53){\makebox(0,0){$C_{g_{32}}$}}

\put(134,53){\makebox(0,0){$R_{g_{33}}$}}
\put(158,53){\makebox(0,0){$C_{g_{33}}$}}

\put(3,37){\makebox(0,0){$g_{34}$}}
\put(63,37){\makebox(0,0){$g_{35}$}}
\put(123,37){\makebox(0,0){$g_{36}$}}

\put(14,3){\makebox(0,0){$R_{g_{34}}$}}
\put(38,3){\makebox(0,0){$C_{g_{34}}$}}

\put(74,3){\makebox(0,0){$R_{g_{35}}$}}
\put(98,3){\makebox(0,0){$C_{g_{35}}$}}

\put(134,3){\makebox(0,0){$R_{g_{36}}$}}
\put(158,3){\makebox(0,0){$C_{g_{36}}$}}

\put(15,80){\line(1,-1){20}}
\put(15,80){\line(2,-1){20}}

\put(15,70){\line(2,1){20}}
\put(15,70){\line(2,-1){20}}

\put(15,60){\line(1,1){20}}
\put(15,60){\line(2,1){20}}

\put(75,80){\line(1,-1){20}}

\put(75,70){\line(1,0){20}}
\put(75,70){\line(2,1){20}}
\put(75,70){\line(2,-1){20}}

\put(75,60){\line(1,1){20}}
\put(75,60){\line(2,1){20}}

\put(135,80){\line(1,-1){20}}
\put(135,80){\line(2,-1){20}}

\put(135,70){\line(1,0){20}}
\put(135,70){\line(2,-1){20}}

\put(135,60){\line(1,1){20}}
\put(135,60){\line(2,1){20}}

\put(15,30){\line(1,-1){20}}

\put(15,20){\line(1,0){20}}
\put(15,20){\line(2,-1){20}}

\put(15,10){\line(1,0){20}}
\put(15,10){\line(1,1){20}}
\put(15,10){\line(2,1){20}}

\put(75,30){\line(1,0){20}}
\put(75,30){\line(1,-1){20}}
\put(75,30){\line(2,-1){20}}


\put(75,10){\line(1,0){20}}
\put(75,10){\line(1,1){20}}
\put(75,10){\line(2,1){20}}

\put(135,30){\line(1,0){20}}
\put(135,30){\line(1,-1){20}}

\put(135,20){\line(2,1){20}}
\put(135,20){\line(2,-1){20}}

\put(135,10){\line(1,0){20}}
\put(135,10){\line(1,1){20}}

\end{picture}
\caption{$n=3$, $k=6$}\label{n3k6}
\end{figure}

For graph $g_{31} \in\mathfrak{G}_{3,6}$ we have:

$$\left[ g_{31} \right] =\left\{ 1,1,1,1,1,1\right\}$$

$$\Psi (g_{31} ) =\langle 0,0,6,0\rangle$$

$$\omega (g_{31} )=\frac{\left[ (3-0)!\right]^0 \left[ (3-1)!\right]^0 }{1!\; 1!\; 1!\; 1!\; 1!\; 1!}  =1$$

For graphs $g_{32} \in \mathfrak{G}_{3,6}$ and $g_{33} \in \mathfrak{G}_{3,6}$ we have:

$$\left[ g_{32} \right] =\left[ g_{33} \right] =\left\{ 1,1,1,1,2\right\}$$

$$\Psi (g_{32} ) =\Psi (g_{33} ) =\langle 0,1,4,1\rangle$$

$$\omega (g_{32} )=\omega (g_{33} )=\frac{\left[ (3-0)!\right]^0 \left[ (3-1)!\right]^1 }{1!\; 1!\; 1!\; 1!\; 2!} = \frac{6^0 \cdot 2^1}{2} =1$$

For graph $g_{34} \in\mathfrak{G}_{3,6}$ we have:

$$\left[ g_{34} \right] =\left\{ 1,1,1,1,1,1\right\}$$

$$\Psi (g_{34} ) =\langle 0,2,2,2\rangle$$

$$\omega (g_{34} )=\frac{\left[ (3-0)!\right]^0 \left[ (3-1)!\right]^2 }{1!\; 1!\; 1!\; 1!\; 1!\; 1!} =\frac{6^0 \cdot 2^2}{1}  =4$$

For graphs $g_{35} \in \mathfrak{G}_{3,6}$ and $g_{36} \in \mathfrak{G}_{3,6}$ we have:

$$\left[ g_{35} \right] =\left[ g_{36} \right] =\left\{ 1,2,3\right\}$$

$$\Psi (g_{35} ) =\Psi (g_{36} ) =\langle 1,0,3,2\rangle$$

$$\omega (g_{35} )=\omega (g_{36} )=\frac{\left[ (3-0)!\right]^1 \left[ (3-1)!\right]^0 }{1!\; 2!\; 3!} = \frac{6^1 \cdot 2^0}{2\cdot 6} =\frac{1}{2}$$

Then for the set $\mathfrak{G}_{3,6}$ we get:

\begin{equation}\label{theta36}
\theta (3,6)=\sum_{g\in \mathfrak{G}_{3,6}} \omega (g)=1+1+1+4+\frac{1}{2} +\frac{1}{2} =8
\end{equation}

\subsubsection{$k=7$}

When $n = 3$ and $k = 7$ the set $\mathfrak{G}_{3,7}= \{ g_{37} ,g_{38} ,g_{39}  \}$ consists of three bipartite graphs, which are shown in Figure \ref{n3k7}.

\begin{figure}[h]
\begin{picture}(170,40)

\multiput(0,0)(60,0){3}{
\put(0,0){\framebox(50,40)}

\multiput(13,21)(24,0){2}{\oval(14,30)}

\put(13,30){\circle*{2}}
\put(37,30){\circle{2}}

\put(13,20){\circle*{2}}
\put(37,20){\circle{2}}

\put(13,10){\circle*{2}}
\put(37,10){\circle{2}}
}

\put(3,37){\makebox(0,0){$g_{37}$}}
\put(63,37){\makebox(0,0){$g_{38}$}}
\put(123,37){\makebox(0,0){$g_{39}$}}

\put(14,3){\makebox(0,0){$R_{g_{37}}$}}
\put(38,3){\makebox(0,0){$C_{g_{37}}$}}
\put(74,3){\makebox(0,0){$R_{g_{38}}$}}
\put(98,3){\makebox(0,0){$C_{g_{38}}$}}
\put(134,3){\makebox(0,0){$R_{g_{39}}$}}
\put(158,3){\makebox(0,0){$C_{g_{39}}$}}

\put(15,30){\line(2,-1){20}}
\put(15,30){\line(1,-1){20}}
\put(15,20){\line(2,1){20}}
\put(15,20){\line(2,-1){20}}
\put(15,10){\line(1,0){20}}
\put(15,10){\line(1,1){20}}
\put(15,10){\line(2,1){20}}

\put(75,30){\line(1,-1){20}}
\put(75,20){\line(2,1){20}}
\put(75,20){\line(1,0){20}}
\put(75,20){\line(2,-1){20}}
\put(75,10){\line(1,0){20}}
\put(75,10){\line(1,1){20}}
\put(75,10){\line(2,1){20}}

\put(135,30){\line(2,-1){20}}
\put(135,30){\line(1,-1){20}}
\put(135,20){\line(1,0){20}}
\put(135,20){\line(2,-1){20}}
\put(135,10){\line(1,0){20}}
\put(135,10){\line(1,1){20}}
\put(135,10){\line(2,1){20}}

\end{picture}
\caption{$n=3$, $k=7$}\label{n3k7}
\end{figure}

For graph $g_{37} \in \mathfrak{G}_{3,7}$ it is true:

$$\left[ g_{37} \right] =\left\{ 1,1,1,1,1,1\right\}$$

$$\Psi (g_{37} ) =\langle 0,0,4,2\rangle$$

$$\omega (g_{37} )=\frac{\left[ (3-0)!\right]^0 \left[ (3-1)!\right]^0 }{1!\; 1!\; 1!\; 1!\; 1!\; 1!} =\frac{6^0 \cdot 2^0}{1}  =1$$

For graphs $g_{38} \in \mathfrak{G}_{3,7}$ and $g_{39} \in \mathfrak{G}_{3,7}$ we get:

$$\left[ g_{38} \right] =\left[ g_{39} \right] =\left\{ 1,1,2,2\right\}$$

$$\Psi (g_{38} ) =\Psi (g_{39} ) =\langle 0,1,2,3\rangle$$

$$\omega (g_{38} )=\omega (g_{39} )=\frac{\left[ (3-0)!\right]^0 \left[ (3-1)!\right]^1 }{1!\; 1!\; 2!\; 2!} = \frac{6^0 \cdot 2^1}{2^2} =\frac{1}{2}$$

Then for the set $\mathfrak{G}_{3,7}$ we get:

\begin{equation}\label{theta37}
\theta (3,7)=\sum_{g\in \mathfrak{G}_{3,7}} \omega (g)=1+\frac{1}{2} +\frac{1}{2} =2
\end{equation}

\subsubsection{$k=8$}

Graph $g_{40}$, which is displayed in Figure \ref{n3k8} is the only bipartite graph belonging to the set $\mathfrak{G}_{3,8}$ in the case $n=3$ and $k=8$.

\begin{figure}[h]
\begin{picture}(50,40)

\put(0,0){\framebox(50,40)}
\put(3,37){\makebox(0,0){$g_{40}$}}

\multiput(13,21)(24,0){2}{\oval(14,30)}
\put(14,3){\makebox(0,0){$R_{g_{40}}$}}
\put(38,3){\makebox(0,0){$C_{g_{40}}$}}

\put(13,30){\circle*{2}}
\put(37,30){\circle{2}}
\put(15,30){\line(1,-1){20}}
\put(15,30){\line(2,-1){20}}

\put(13,20){\circle*{2}}
\put(37,20){\circle{2}}
\put(15,20){\line(1,0){20}}
\put(15,20){\line(2,1){20}}
\put(15,20){\line(2,-1){20}}

\put(13,10){\circle*{2}}
\put(37,10){\circle{2}}
\put(15,10){\line(1,0){20}}
\put(15,10){\line(1,1){20}}
\put(15,10){\line(2,1){20}}

\end{picture}
\caption{$n=3$, $k=8$}\label{n3k8}
\end{figure}

For graph $g_{40} \in \mathfrak{G}_{3,8}$  it is true:

$$\left[ g_{40} \right] =\left\{ 1,1,2,2\right\}$$

$$\Psi (g_{40} ) =\langle 0,0,2,4\rangle$$

$$\omega (g_{40} )=\frac{\left[ (3-0)!\right]^0 \left[ (3-1)!\right]^0 }{1!\; 1!\; 2!\; 2!} =\frac{6^0 \cdot 2^0}{2^2}  =\frac{1}{4}$$

Therefore:

\begin{equation}\label{theta38}
\theta (3,8)=\sum_{g\in \mathfrak{G}_{3,8}} \omega (g)=\frac{1}{4}
\end{equation}

\subsubsection{$k=9$}

When $n = 3$ and $k = 9$ there is only one graph and this is the complete bipartite graph $g_{41}$ which is shown in Figure \ref{n3k9}.

\begin{figure}[h]
\begin{picture}(50,40)

\put(0,0){\framebox(50,40)}
\put(3,37){\makebox(0,0){$g_{41}$}}

\multiput(13,21)(24,0){2}{\oval(14,30)}
\put(14,3){\makebox(0,0){$R_{g_{41}}$}}
\put(38,3){\makebox(0,0){$C_{g_{41}}$}}

\put(13,30){\circle*{2}}
\put(37,30){\circle{2}}
\put(15,30){\line(1,0){20}}
\put(15,30){\line(1,-1){20}}
\put(15,30){\line(2,-1){20}}

\put(13,20){\circle*{2}}
\put(37,20){\circle{2}}
\put(15,20){\line(1,0){20}}
\put(15,20){\line(2,1){20}}
\put(15,20){\line(2,-1){20}}

\put(13,10){\circle*{2}}
\put(37,10){\circle{2}}
\put(15,10){\line(1,0){20}}
\put(15,10){\line(1,1){20}}
\put(15,10){\line(2,1){20}}

\end{picture}
\caption{$n=3$, $k=9$}\label{n3k9}
\end{figure}

For graph $g_{41}$ is true:

$$\left[ g_{41} \right] =\left\{ 3,3\right\}$$

$$\Psi (g_{41} ) =\langle 0,0,0,6\rangle$$

$$\omega (g_{41} )=\frac{\left[ (3-0)!\right]^0 \left[ (3-1)!\right]^0 }{3!\; 3!} =\frac{6^0 \cdot 2^0}{6^2}  =\frac{1}{36}$$

Therefore

\begin{equation}\label{theta39}
\theta (3,9)=\sum_{g\in \mathfrak{G}_{3,9}} \omega (g)=\frac{1}{36}
\end{equation}

Having in mind the formula (\ref{endstar}) and formulas  (\ref{theta31}) $\div$ (\ref{theta39}) for the number $D_9$ of all ordered pairs disjoint S-permutation  matrices in $ n =  3$ we finally get:

\begin{equation}\label{D_9}
D_9 = (3!)^{12} + (3!)^{8} \left[ \sum_{k=1}^9 (-1)^k \theta (n,k) \right] =
\end{equation}
$$= 2\; 176\; 782\; 336 + 1\;679\; 616  \left(-1296+1008-352+125-37+8-2+\frac{1}{4} - \frac{1}{36}\right) = $$
$$=1\; 260\; 085\; 248 . $$

The number $d_9$ of all non-ordered pairs disjoint matrices from $\Sigma_9$  is equal to

\begin{equation}
d_9 =\frac{1}{2} D_9 =630\; 042\; 624
\end{equation}

\subsection{On a combinatorial problem of graph theory related to the number of  Sudoku matrices }

\begin{problem}\label{prbl1}
Let $n\ge 2$ is a natural number and let $G$ be a simple graph having $(n!)^{2n}$ vertices. Let each vertex of $G$ be identified with an element of the  set $\Sigma_{n^2}$ of all $n^2 \times n^2$ S-permutation matrices. Two vertices are connected by an edge if and only if the corresponding matrices are disjoint. The problem is to find the number of all complete subgraphs of $G$ having $n^2$ vertices:
\end{problem}

Note that the number of edges in graph $G$ is equal to $d_{n^2}$ and can be calculated using formula (\ref{f2}) and formula (\ref{f3})  (respectively formulas (\ref{star2}), (\ref{thet}), (\ref{endstar})  and (\ref{f3})).

Denote by $z_n$  the solution of the Problem \ref{prbl1} and let $\sigma_n$ is the number of all $n^2 \times n^2$ Sudoku matrices. Then according to Proposition \ref{disj} and the method of construction of the graph
 $G$, it follows that the next equality is valid:
\begin{equation}\label{z_n}
z_n =\frac{\sigma_n}{(n^2 )!}
\end{equation}

We do not know a general formula for finding the number of all $n^2 \times n^2$ Sudoku matrices for each natural number $n\ge 2$ and we consider that this is an open combinatorial problem. Only some special cases are known. For example in $n=2$ it is known that $\sigma_2 =288$ \cite{yorkost}. Then according to formula (\ref{z_n}) we get:

$$z_2 =\frac{\sigma_2}{4!} =\frac{288}{24} =12$$

In  \cite{Felgenhauer} it has been shown that in $n=3$ there are exactly
$$\sigma_3 = 6\; 670\; 903\; 752\; 021\; 072\; 936\; 960 =$$
$$= 9! \times 72^2 \times 2^7 \times 27\; 704\; 267\; 971 =$$
$$2^{20} \times  3^8 \times  5^1 \times  7^1 \times 27\; 704\; 267\; 971^1  \sim 6.671\times 10^{21}$$
number of Sudoku matrices. Then according to formula (\ref{z_n}) we get:

$$z_3 =\frac{\sigma_3}{9!} =\frac{6\; 670\; 903\; 752\; 021\; 072\; 936\; 960}{362\; 880} =18\; 383\; 222\; 420\; 692\; 992$$

\bibliographystyle{plain}
\bibliography{bigraphs}

\begin{thebibliography}{1}

\bibitem{Bailey}
R.A. Bailey, P.J. Cameron, and R.~Connelly.
\newblock Sudoku, gerechte designs, resolutions, affine space, spreads, reguli,
  and hamming codes.
\newblock {\em Amer. Math. Monthly}, (115):383--404, 2008.

\bibitem{dahl}
G.~Dahl.
\newblock Permutation matrices related to {Sudoku}.
\newblock {\em Linear Algebra and its Applications}, (430):2457--2463, 2009.

\bibitem{diestel}
R.~Diestel.
\newblock {\em Graph Theory}.
\newblock Springer-Verlag Heidelberg, New York, 1997, 2000, 2006.

\bibitem{Felgenhauer}
B.~Felgenhauer and F.~Jarvis.
\newblock Enumerating possible {Sudoku} grids.
\newblock 2005.
\newblock http://www.afjarvis.staff.shef.ac.uk/sudoku/sudoku.pdf.

\bibitem{Fontana}
R.~Fontana.
\newblock Fraction of permutations - an application to sudoku.
\newblock {\em Journal of Statistical Planning and Inference},
  (141):3697--3704, 2011.

\bibitem{harary}
F.~Harary.
\newblock {\em Graph Theory}.
\newblock Addison-Wesley, Massachusetts, 1998.

\bibitem{yorkost}
H.~Kostadinova and K.~Yordzhev.
\newblock An entertaining example for the usage of bitwise operations in
  programming.
\newblock In {\em FMNS-2011}, volume~1, pages 159--168, Blagoevgrad, Bulgaria,
  2011. SWU.

\bibitem{KYbigraphs}
K.~Yordzhev.
\newblock On the number of disjoint pairs of s-permutation matrices.
\newblock 2012.
\newblock arXiv:1211.1628.

\end{thebibliography}

\end{document}